\documentclass[12pt,leqno]{article}
\usepackage{amsmath,amssymb,amsfonts}
\usepackage{eucal}

\newcommand\comment[1]{}                

\newcommand\was[1]{\textbf{OLD} \texttt{[#1]}}
\renewcommand\was[1]{}
\newcommand\now[1]{\textbf{NEW} [#1]}
\renewcommand\now[1]{#1}

\numberwithin{equation}{section}

\newtheorem{lem}{Lemma}[section]
\newtheorem{cor}[lem]{Corollary}
\newtheorem{prop}[lem]{Proposition}

\newtheorem{theorem}[lem]{Theorem}
\newtheorem{conjecture}[lem]{Conjecture}
\newtheorem{corollary}[lem]{Corollary}
\newtheorem{lemma}[lem]{Lemma}

\newenvironment{proof}[1][Proof]{\textbf{#1.} }{\ \rule{0.5em}{0.5em}\smallskip}
\newtheorem{example}{Example}[section]

\renewcommand{\phi}{\varphi}                 
\renewcommand{\epsilon}{\varepsilon}
\renewcommand\tilde{\widetilde}

\newcommand\Aut{\operatorname{Aut}}

\newcommand\rank{\operatorname{rank}}
\renewcommand\Im{\operatorname{Im}}           

\newcommand\Ker{\operatorname{Ker}}

\newcommand\lra{\longrightarrow}
\newcommand\Sat{\operatorname{Sat}}
\newcommand\SatPair{\operatorname{SatPair}}
\newcommand\Bal{\operatorname{Bal}}
\newcommand\Star{\operatorname{Star}}
\newcommand\Lift{\operatorname{Lift}}
\newcommand\Rec{\operatorname{Rec}}
\newcommand\Sk{\operatorname{Sk}}

\newcommand\G{\Gamma}
\newcommand\GG{(\G,g,G)}

\renewcommand\b{\mathbf{b}}
\renewcommand\c{\mathbf{c}}
\newcommand\e{\mathbf{e}}
\newcommand\g{\mathbf{g}}
\newcommand\h{\mathbf{h}}
\newcommand\n{\mathbf{n}}

\newcommand\w{\mathbf{w}}
\newcommand\x{\mathbf{x}}
\newcommand\y{\mathbf{y}}
\newcommand\z{\mathbf{z}}

\newcommand\cG{\mathcal{G}}
\newcommand\cH{\mathcal{H}}
\newcommand\cK{\mathcal{K}}
\newcommand\cM{\mathcal{M}}

\newcommand\cR{\mathcal{R}}

\newcommand\R{\mathbb{R}}    
\newcommand\RP{\mathbb{RP}} 
\newcommand\Z{\mathbb{Z}}    
\newcommand\bbR{\mathbb{R}}
\newcommand\bbZ{\mathbb{Z}}

\begin{document}

\title{Criteria for Balance in Abelian Gain Graphs, \\
    with Applications to Piecewise-Linear Geometry\footnote{We thank Keith Dennis and  Kenneth S.\ Brown for very helpful discussions.}
    }
\author{Konstantin Rybnikov\footnote{Department of Mathematical Sciences, University of Massachusetts at Lowell, Lowell, MA 01854, U.S.A.  Part of this research was conducted at Cornell University, Ithaca, NY 14853-4201, U.S.A. Research partially supported by the Mathematical Sciences Research Institute in Berkeley.}  \\
    and Thomas Zaslavsky\footnote{Department of Mathematical Sciences, Binghamton University, Binghamton, NY 13902-6000, U.S.A.  Research partially supported by National Science Foundation grant DMS-0070729.}
    }
\date{\today}

\maketitle

\emph{Short title}: Balance in Abelian Gain Graphs

\abstract{Consider a gain graph with abelian gain group having no odd
torsion.  If there is a basis of the graph's binary cycle space each of
whose members can be lifted to a closed walk whose gain is the identity,
then the gain graph is balanced, provided that the graph is finite or the
group has no nontrivial infinitely 2-divisible elements.  We apply this
theorem to deduce a result on the projective geometry of
piecewise-linear realizations of cell-decompositions of manifolds.}

\bigskip

\emph{Keywords}: abelian group, gain graph, balanced gain graph, binary
cycle, Binary Cycle Test, cell-complex, cell-decomposition of manifold,
dual graph, facet graph,  homology manifold, polyhedral manifold,
piecewise-linear surface, multivariate spline, convex surface,  Maxwell
lifting, reciprocal graph, Maxwell--Cremona correspondence.

\emph{Mathematics Subject Classifications (2000)}:
\emph{Primary} 05C22, 52B70; \emph{Secondary} 05C38, 52C25, 52C22.
\comment{05C22 Gain [etc.] graphs,
52B70 Geometry of Polyhedral Manifolds,
05C38 Paths and cycles,
52C25 Rigidity and flexibility of structures,
52C22 Tilings in $n$ dimensions}

\bigskip\hrule\bigskip

This paper does NOT have a ``corresponding author'' or ``senior author''.
Both authors are able to answer correspondence from readers.

\bigskip

\bigskip\hrule\bigskip


\section{Introduction}

A \emph{gain graph} $\GG$ consists of a graph $\G=(V,E)$, a group $G$,
and a homomorphism $g$ from the free group $F_E$ on the edges of $\G$ to
$G$.  We call $g$ the \emph{gain map} and $G$ the \emph{gain
group}.\footnote{See, e.g., Zaslavsky (1989), Section 5, for the basic
theory of gain graphs.  Gain graphs have been called \emph{voltage
graphs} in graph embedding theory---e.g., in Gross \& Tucker (1987); we
eschew the term ``voltage'' because gains do not have to obey Kirchhoff's
voltage law.}
(One thinks of the edges of $G$ as oriented in an arbitrary but fixed
way).
Gain graphs have appeared in pure mathematics, physics, operations
research, psychology, and discrete geometry, among many other areas (see
Zaslavsky's annotated bibliography (1998) for a survey).

The simplest gain graphs are those in which every simple closed walk on $\Gamma$ lies in the kernel of $g$; such gain graphs are called \emph{balanced}.
If $\Gamma$ is connected, $\GG$ is balanced if and only if $\pi(\G) \subseteq \Ker g$.
Since for general gain graphs one can treat each component separately, we shall assume throughout that $G$ is connected.
A basic problem is how to tell whether or not a gain graph is balanced.
We investigate this problem when $G$ is abelian;  then we call $\GG$ an \emph{abelian gain graph}.
We propose here a criterion for balance in abelian gain graphs that is broader than previous ones, but which has the drawback that it is not valid for all groups.  Our topics are, first, the problem of describing the abelian gain groups for which it is valid and, second, an application of the criterion to a problem of polyhedral geometry.

To state the criterion we must first define a ``circle'' and a ``cycle'' in a graph.  A \emph{circle} is the edge set of a nontrivial simple closed walk (a walk in which no vertex or edge is repeated, except that the initial and final vertices are the same).
There are two kinds of cycles of concern to us:\footnote{The graph-theory literature contains at least five incompatible definitions of the term ``cycle''.}
they are \emph{integral cycles}, i.e., elements of $Z_1(\G) = \Ker \partial$ where $\partial$ is the standard 1-dimensional boundary operator of cellular homology of $\G$, and \emph{binary cycles}, i.e., elements of $Z_1(\G;\Z_2)$.
There are natural epimorphisms
$$
\pi(\G) \overset{a}{\lra} Z_1(\G) \overset{[\ ]}{\lra} Z_1(\G;\Z_2)
$$
that consist, respectively, of abelianization and reduction of coefficients modulo 2.
Suppose that $b$ is a binary cycle: a \emph{cyclic orientation} of $b$ is any closed walk $\tilde b$ for which $[\tilde b^{a}] = b$.  If $B$ is a set of binary cycles, a \emph{cyclic orientation of $B$} is any set $\tilde B = \{ \tilde b : b \in B \}$ of cyclic orientations of the members of $B$.  This definition is very broad: a binary cycle has many cyclic orientations; also, edges may appear in $\tilde b$ that are not in the support of $b$, provided that each such edge appears an even number of times in $\tilde b$.  Thus, a binary cycle with disconnected support can have a cyclic orientation, which is always connected. In general, there is no cyclic orientation of a cycle that can reasonably be regarded as canonical, except in the case of a binary cycle whose support is a circle.

Now, our test for balance:

\medskip
\noindent \textbf{Definition. }
Let $\GG$ be a gain graph and $B$
be a basis of $Z_1(\G;\Z_2)$.  We say that $B$ \emph{passes the Binary
Cycle Test} if it has a cyclic orientation $\tilde B$ such that all
elements of $\tilde B$ have gain 1. We say the Binary Cycle Test is
\emph{valid for $\GG$} if the existence of a basis $B$ that passes the
Binary Cycle Test implies that $\GG$ is balanced. We say the Binary Cycle
Test is \emph{valid} for a family of graphs $\cG$ and a family of
groups $\cH$ if it is valid for every gain graph $\GG$ with
$\G \in \cG$ and $G \in \cH$.
\medskip

Our main results are that the Binary Cycle Test is valid for any abelian
gain group that has no odd torsion if the graph is finite (Theorem
\ref{T:finite-abelian}) or if the gain group has no nontrivial infinitely
2-divisible elements (Theorem \ref{T:abelian}).  We also show an infinite
graph that fails the Binary Cycle Test for all gain groups that have an
infinitely 2-divisible element other than the identity.  (In Rybnikov \&
Zaslavsky (20xx) we prove that for any gain group with non-trivial
elements of odd order there is a graph for which the Binary Cycle Test
does not work.  Focussing on properties of the gain graph that are due to
the gain group rather than the graph is how the present paper differs from
that one.)

In most applications to physics, psychology, or operations research, the
gain graphs are finite and have $\Z_2$ or $\R^*$ as the gain group. Groups
$\R^*$ or $\R^d$ tend to appear in discrete-geometry applications such as
Voronoi's generatrix construction and the Maxwell--Cremona correspondence.
The Binary Cycle Test applies to such situations, since all these groups
are abelian and do not have odd torsion. However, if the gain graph is
infinite, our results do not apply when $G = \R^*$ or $\R^d$, since these
groups are 2-divisible.  We know that for certain particular binary cycle
bases their passing the Binary Cycle Test does imply balance even then,
but they have to be carefully selected, as Example \ref{X:infinite} shows.

The reason we need the Binary Cycle Test rather than a simpler test for
balance of $\GG$, such as the existence of a set of balanced circles that
form a fundamental system of circles (i.e., the set of circles obtained
from a fixed spanning tree by adjoining each nontree edge), is that, in
applications, often a basis of the binary cycle space is given to us or is
in some way naturally associated with the problem, so we are not free to
choose the basis, or the geometry may yield information about the gains
only of certain cycles that are not a fundamental system or otherwise
tractable.

\medskip
In Euclidean space, the sphere, or hyperbolic space, gain graphs are a natural setting
for problems of existence and classification of tilings with prescribed restrictions on
the group or the tiles, and of $d+1$-dimensional piecewise-linear realization of
$d$-manifolds with given projection.  Among important examples are Coxeter's theorem
reducing the classification of reflection groups to that of Coxeter simplices (1934),
Alexandrov's space-filling theorem (1954), Voronoi's (1908) construction of a
generatrix (\emph{generatice} in the French original), the Maxwell--Cremona
correspondence and its generalizations (see Crapo \& Whiteley, 1993; Rybnikov, 1999;
Erdahl, Rybnikov, \& Ryshkov, 2001), colorings of graphs and polyhedra (Ryshkov \&
Rybnikov, 1997), etc. In Section \ref{applied} we demonstrate the appropriateness of
gain graphs by applying the Binary Cycle Test to generalize a theorem of Voronoi
(1908), that any simple tiling of $\R^d$ is the projection of a (convex) PL-surface
which, given the tiling, is uniquely determined by $d+2$ free real parameters. (A
spline theorist would say that the dimension of the space of $C_1^0$-splines over the
tiling is $d+2$.) Such a surface is called a \emph{(convex) lifting of the tiling to}
$\R^{d+1}$.  It is known that the notion of lifting and the space of liftings can be
defined for any PL-realization of a homology $d$-manifold in $\R^d$. Our generalization
consists in showing that the dimension of the space of liftings of any non-degenerate
PL-realization in $\R^{d}$ of a homology $d$-manifold $\cM$ with $H_1(\cM;\Z_2)=0$ is
$d+2$, provided that each $(d-3)$-cell of $\cM$ is incident to exactly four $d$-cells.
This application was chosen because in the course of the proofs gain graphs appear in a
number of different geometric guises and we employ the Binary Cycle Test several times.

\comment{Consider removing this if we are asked to further shorten.}
Another application of the Binary Cycle Test is to Voronoi's conjecture on
parallelohedra, which states that any parallelohedron is an affine copy of
the Dirichlet domain of a lattice. Deza and Grishukhin (2003) recently
gave a number of combinatorial reformulations of Voronoi's conjecture that
are essentially gain-graphical. The proof of the Voronoi conjecture can be
reduced to a number of questions about a certain type of finite gain graph
with gain group $\bbR^*$, called the \emph{Venkov graph}.  (See Ordine
(2002), Rybnikov (2002a,b), and Deza \& Grishukhin (2003) for definitions
and examples.) We have noticed that the Binary Cycle Test simplifies
attacking these questions, because in the case of Venkov graphs it is much
easier to look for a basis of the binary cycle space rather than of the
integral cycle space. Very recently Ordine (2004) made a major advance in
Voronoi's conjecture by proving an important special case, where the dual
figures for the stars of $(d-3)$-cells of the tiling are restricted to
tetrahedra, octahedra, and pyramids. The Binary Cycle Test can be applied
to greatly simplify his proof; we hope to explain this elsewhere.


\section{Preliminaries}

Our graphs may be infinite.  Loops (edges with two equal endpoints) and
multiple edges are allowed.  We assume that $\G$
is connected and that it has a fixed (but arbitrary) orientation.

A circle $C$ can be written as the edge set of a simple closed walk in a way that is unique up to choice of initial vertex and direction.  Then it depends only on $C$, not on the choice of walk, whether the walk is in $\Ker g$.  If it is, we say $C$ is \emph{balanced}.

If $\tilde{w}$ is a closed walk on $\Gamma$, denote by $\w$ the corresponding element ${\tilde w}^{a}$ of $Z_{1}(\G)$.  If $\w \in Z_1(\G)$, denote by $[\w]$ the corresponding element of $Z_1(\G;\Z_2)$.

We denote by $\pi(\G)$ the fundamental group of $\G$.
The \emph{essential gain group} of a gain graph is $g(\pi(\G))$.
It is well defined up to conjugation in $G$.
Evidently, $\GG$ is balanced if and only if the essential gain group is trivial.

Naturally, $g$ induces a homomorphism from $\pi(\G)$ to $G$ and from $\Z E$, the free abelian group on $E$, to $G/G'$, where, as usual, $G'$ denotes the commutator subgroup.  Hence, $g$ induces a homomorphism from $Z_1(\G)$, the abelianization of $\pi(\G)$, to $G/G'$.

A \emph{cyclic orientation}, or \emph{lifting to the fundamental group}, of an element $b$ of $Z_{1}(\G;\Z_2)$ is any element $\tilde{b}$ of $\pi(\G)$ whose abelianization modulo 2 is $b$.
A circle has natural cyclic orientations obtained by taking an element of $\pi(\G)$ corresponding to a simple walk around the circle in either direction from any starting point.
A \emph{lifting of $b$ to the cycle space} is the image $\b$ in $Z_1(\G)$ of a cyclic orientation of $b$.
If $B$ is a basis of $Z_1(\G;\Z_2),$ then denote by $\tilde{B}$ a cyclic orientation of $B$.
As one can see from our notation, if $b$ denotes an element of the binary cycle space, then $\b$ denotes one of the possible liftings of $b$ to the cycle space, and $\tilde{b}$ is one of the liftings of $b$ to the fundamental group. If $\w \in Z_1(\G)$, denote by $[\w]$ the corresponding element of $Z_1(\G;\Z_2)$.

If $H \leq \pi(\G)$, then $H^a$ denotes $H/(\pi(\G)'\cap H)$, the abelianization of $H$ with respect to $\pi(\G)$; that is, $H^a$ is the image of $H$ in $Z_1(\G)$.
If $L, M \leq A$, we denote by $L/M$ the image of $L$ in $A/M$; in particular, $L/2A$ is $L$ reduced modulo $2$.


\section{The Test for Abelian Groups}

\begin{lemma} \label{L:free-abelian}
Let $A$ be a free abelian group, and let $h$ be a homomorphism from $A$ to $G$,
an abelian group without odd torsion.  Let $L$ be a subgroup of $A$ such
that $h(L)=0$ and $L/2A = A/2A$.  Suppose $K$ is a subgroup of $A$
such that $K/(L\cap K)$ does not have elements of infinite order.  Then $h(K)=0$.
\end{lemma}

\begin{proof}
Suppose there is $\y\in K$ such that $h(\y)\neq 0$.  Since
$h(L)=0$, $\y\notin L$.  There is a minimal $l>1$ such that $l\y\in L$.
Let $l=2^{e}a$, where $a$ is odd. Denote $a\y$ by $\x$. Notice $e>0$, since otherwise $h(\x)$ would have odd order in $G$.

Now, $h(a\y)\neq 0$, since otherwise $h(\y)$ would have odd order in $G$.   Then $2^{e}$ is the order of $\x$ in $K/(L\cap K)$ and
$e = \min \{k \mid 2^{k}\x\in L\}$.  Let $\{\b_{i}\}$ be a basis of $L$, so that
$2^{e}\x = \sum l_{i}\b_{i}$.  Since $L/2A=A/2A$, $\{[\b_{i}]\}_i$ generates
$A/2A$.  Thus $[\x] = \sum \lambda _{i}[\b_{i}]$, i.e.,
$\x = \sum \lambda_{i}\b_{i}+2\z$ for some $\z\in A$.
We have $2^{e}\x = \sum \lambda _{i}2^{e}\b_{i}+2^{e+1}\z$.  Thus
$\sum (l_{i}-\lambda _{i}2^{e})\b_{i}\in 2A$ and all $l_{i}$ are divisible by $2$, i.e., $l_{i}=2m_{i}$, $m_{i}\in \Z$.  Since $e>0$ we have $2^{e-1}\x = \sum m_{i}\b_{i}\in L$.  This contradicts the choice of $e$ as $\min \{k \mid 2^{k}\x\in L\}$.  Therefore $h(K)=0.$
\end{proof}

\begin{lemma} \label{L:torsion-test-for-gains}
Let $\GG$ be a gain graph, where $G$ is an abelian group without odd torsion,
and let $K$ be a subgroup of $Z_{1}(\Gamma)$.
If a basis $B$ of $Z_{1}(\Gamma;\Z_2)$
has some cyclic orientation $\tilde{B}$ that lies in the kernel of $g$,
and if $K/(\langle\tilde{B}\rangle^{a} \cap K)$ does not have elements of infinite
order, then any lifting of $K$ to $\pi(\Gamma)$ lies in $\Ker g$.
\end{lemma}

\begin{proof}
Denote $\langle \tilde{B}\rangle ^{a}$ by $L$.  The homomorphism $g$ induces, in
a natural way, a homomorphism $h$ from $Z_{1}(\Gamma)$ to $G$. Obviously,
$h(L)=0$ and $L/2Z_{1}(\Gamma) = Z_{1}(\Gamma)/2Z_{1}(\Gamma) = Z_{1}(\Gamma;\Z_2).$  We can now apply Lemma \ref{L:free-abelian} and
conclude that $K$ lies in the kernel of $h$.  Since $G$ is abelian, any
lifting of $K$ to $\pi(\Gamma)$ also lies in $\Ker g$.
\end{proof}

\begin{theorem} \label{T:finite-abelian}
Suppose $\Gamma$ is a \emph{finite} graph and $G$ is an abelian group
without odd torsion. If all elements of a basis $B$ of $Z_{1}(\Gamma;\Z_2)$
have cyclic orientations that lie in $\Ker g$,
then $\GG$ is balanced.
\end{theorem}

\begin{proof}
Since $\Gamma$ is finite,
$\rank\langle \tilde{B}\rangle ^{a} = \rank Z_{1}(\Gamma)$,
so all elements of $Z_{1}(\Gamma)/\langle \tilde{B}\rangle^{a}$ have finite order.
Now apply Lemma \ref{L:torsion-test-for-gains} with $K=Z_{1}(\Gamma)$.
\end{proof}

\begin{lemma} \label{L:no-infinite-order}
Let $G$ be an abelian group without odd torsion. If some basis $B$ of
$Z_{1}(\Gamma;\Z_2)$ has a cyclic orientation $\tilde{B}$ that lies
in $\Ker g$, then $g(\pi(\Gamma))$ is 2-divisible.
\end{lemma}

\begin{proof}
The elements of $Z_{1}(\Gamma)$ that have finite order in
$Z_1(\Gamma)/\langle\tilde{B}\rangle^{a}$ form a subgroup in $Z_1(\Gamma)$, which we
denote by $F$. Since $F/\langle \tilde{B}\rangle ^{a}$ does not have elements of
infinite order, Lemma \ref{L:torsion-test-for-gains} with $K = F$ implies that
$g(F)=0$.  Thus if $\y-\x \in F$, $g(\y-\x)=0.$  Therefore, $g$ induces an epimorphism
$h$ from $Z_{1}(\Gamma)/F$ onto $g(\pi(\Gamma))$.

Let $\y\in Z_{1}(\Gamma).$  Since $\langle \tilde{B}\rangle ^{a} \subset F$ and $[\langle \tilde{B}\rangle ^{a}] = Z_1(\G,\Z_2)$,
there is $\x \in F$ such that $\y-\x = 2\z$ for some $\z\in Z_{1}(\Gamma)$.
Thus $\y\in 2\z+F,$ and $\y+F = 2\z+F= 2(\z+F)$ as posets.  Therefore all elements
of $Z_{1}(\Gamma)/F$ are 2-divisible.

Let $a \in g(\pi(\Gamma))$.  Since the homomorphism $h$ from
$Z_{1}(\Gamma)/F$ is onto $g(\pi(\Gamma))$, there is $\y\in Z_{1}(\Gamma)/F$ such that $h(\y)=a$.  Since $\y$ is 2-divisible, $a$ is also 2-divisible: $a=h(\y)=h(2\z)=2h(\z).$
\end{proof}

\begin{theorem} \label{T:abelian}
Let $G$ be an abelian group, without odd torsion, which has no
infinitely 2-divisible elements other than $0$.  If some basis $B$ of $Z_{1}(\Gamma;\Z_2)$ has a cyclic orientation $\tilde{B}$ that lies in
$\Ker g$, then $\GG$ is balanced.
\end{theorem}

\begin{proof}
If $Z_{1}(\Gamma)/\langle \tilde{B}\rangle ^{a}$ does not have elements of
infinite order, then, an application of Lemma \ref{L:torsion-test-for-gains}
with $K = Z_{1}(\Gamma)$ implies that $\GG$ is balanced.
If $Z_{1}(\Gamma)/\langle \tilde{B}\rangle ^{a}$ has
elements of infinite order, then by Lemma \ref{L:no-infinite-order}
$g(\pi(\Gamma))$ is 2-divisible, which is only possible if $g(\pi(\Gamma))=0$. Thus $\GG$ is balanced.
\end{proof}

\begin{corollary} \label{T:abelian-finite-free}
Suppose $\Gamma$ is a \emph{finite} graph and $G$ is a free group. If all
elements of a basis $B$ of $Z_1(\Gamma;\Z_2)$ have cyclic
orientations that lie in $\Ker g$, then $\GG$ is balanced.
\end{corollary}

\begin{proof}
$g$ induces a homomorphism $\overline{g}$ from $\pi(\Gamma)$ into
$G/G^{\prime}.$  $G/G^{\prime}$ is a free abelian group and, therefore, does
not have torsion or infinitely 2-divisible elements.  By
Theorem \ref{T:abelian}, $\Im\overline{g}=0.$
But this can only happen if $\Im g=1$, because $\Ker g$ is a free group whose abelianization is $\Ker\overline{g}$.
\end{proof}

\begin{cor} \label{T:nonabelian}
Let $g$ be a group such that $G/G'$ has no torsion or has only 2-torsion.  Let $\GG$\ be a gain graph whose gain group is $G$.  If there is a basis $B$ of the binary cycle space $Z_1(\G;\Z_2)$ such that for each $b\in B$ there is a cyclic orientation $\tilde b$ for which $g(\tilde b) \in G'$, then every circle has gain in $G'$.
\end{cor}


\section{Counterexamples}  \label{counter}

In discrete geometry the case of an infinite gain graph with an abelian
2-divisible group of gains, such as $\R^n$  or $\R^*$ is one of the most
typical ones (see Ryshkov \& Rybnikov (1997), Rybnikov (1999)).
In such applications the cycle basis for which we can verify the balance
property  normally consists of \emph{circles}. Might the \emph{Circle
Test}, i.e., the Binary Cycle Test where all elements of the binary basis
are known to be \emph{circles}, work for infinite graphs with 2-divisible
groups of gains?  Unfortunately, not.  We construct an infinite gain graph
with a 2-divisible, torsion-free group of gains, unbalanced yet having a
binary circle basis all of whose elements are balanced.

\emph{A more detailed description of examples considered, or only mentioned, in this
section is given in Appendix, following the paper.}

\begin{example} \label{X:infinite}{\rm
The graph is shown in Figure 1.
We use the following notation for elements of $Z_1(\G)$:
$$
\e_k = D_kL_kD_{k+1}U_{k}D_k, \quad \g_k = L_kL_{k+1}D_{k+1}L_k, \quad
\h_k = U_kD_{k+1}U_{k+1}U_k.
$$
Then $\{\e_k, \g_{k}, \h_{k} \}_{k \ge 1}$ is an integral circle basis of
$\G$.  As usual, we denote an element of $Z_1(\G;\Z_2)$ corresponding to
$\x \in Z_1(\G)$ by $[ \x ]$.  Let
$$
\mathbf{B} =
\{\e_{k}+\g_{k}+\h_{k}\}_{k \ge 1} \cup \{\g_k+\e_{k+1} \}_{k \ge 1}
\cup \{\h_k+\e_{k+1}\}_{k \ge 1}
$$
and let $B = \{[\b] \mid \b \in \mathbf{B}\}$.  $B$ is obviously a basis of
$Z_1(\G;\Z_2)$, because $[\e_k] =
[\e_{k}+\g_{k}+\h_{k}]+[\g_{k}+\e_{k+1}]+[\h_{k}+\e_{k+1}]$. Let $H$ denote the
subgroup of $Z_1(\G)$ (hence of $\Z E$) generated by $\mathbf{B}$, and let $\chi: \Z E
\to \Z E/H$ be the quotient map. Now, let us set the gain group to be $\Z E/H$ and the
gain map $g$ to be $\chi a$, where $a$ is the abelianization mapping from $F_E$  to $\Z
E$. Then the essential gain group $g(\pi(\G))=\chi(Z_1(\G))=Z_1(\G)/H$. By
construction, $B$ is a basis of the binary cycle space that consists of balanced
circles.  If $H \neq Z_1(\G)$, then the constructed gain graph is not balanced, because
$Z_1(\G)/H = \langle \chi(\e_k), \chi(\g_{k}), \chi(\h_{k}) \mid k \ge 1 \rangle =
g(\pi(\G))$. But in fact, $H \neq Z_1(\G)$, because there is no finite linear
combination of elements of $H$ that is equal to $\e_1$.

\begin{figure}
\vbox to 3.3 truein{}
\begin{center}
\includegraphics{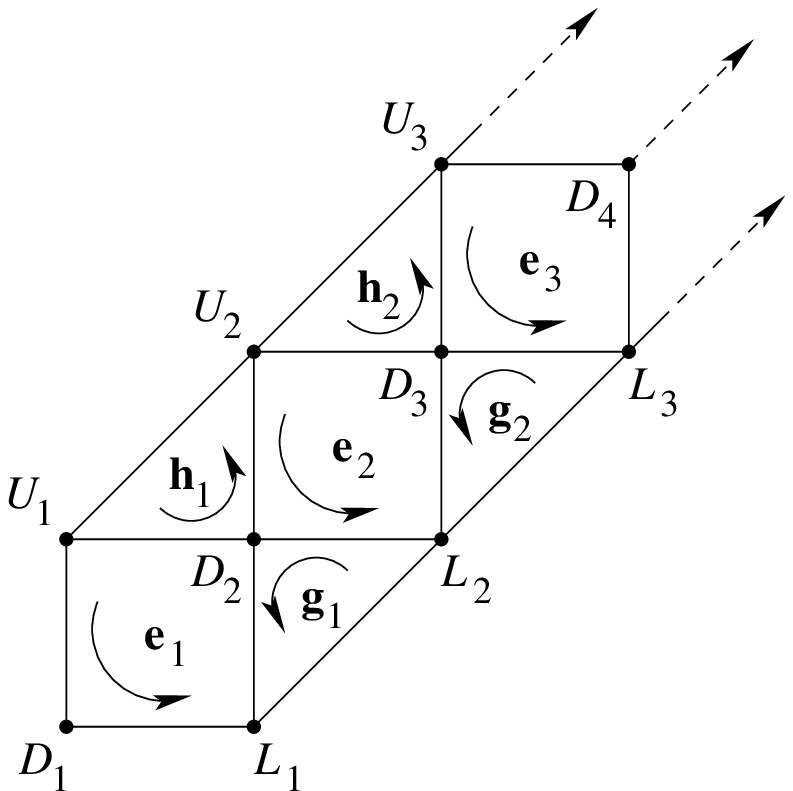}
\end{center}
\caption{Counterexample to the Circle Test with infinitely 2-divisible
elements.}   \label{F:egh-graph}
\end{figure}

Denote by $\Z[\frac{1}{2}]$ the ring of dyadic rationals, i.e., the ring
of rational numbers whose denominators are powers of 2.
Its additive group $\Z[\frac{1}{2}]^+$ is the unique smallest group that
contains an infinitely 2-divisible element of infinite order.
One can show that $Z_1(\G)/H \cong \Z[\frac{1}{2}]^+$ by mapping
$\e_i+H \mapsto d_i$ where $\{d_i\}_{i=1}^\infty$ generates
$\Z[\frac{1}{2}]^+$ and satisfies $d_i = 2 d_{i+1}$.
Thus, the essential gain group $g(\pi(\G))$ is 2-divisible.
}
\end{example}

Notice that a circle $\e_{i}+\g_{i}+\h_{i}$ has 5 edges.

\begin{conjecture}
Let $G$ be a torsion-free abelian group, possibly 2-divisible, and
suppose $Z_1(\G)$ has a basis that consists of balanced circles with at
most $k$ edges. Then $\GG$ is balanced if $k=3$, but it need not be
balanced if $k=4$.
\end{conjecture}

There are also finite counterexamples, but the gain group necessarily has an element of
odd order. Two such examples, based respectively on the wheel and on an even circle
with all edges doubled, are in Rybnikov \& Zaslavsky (20xx), Theorem 6.16 (see Appendix
for one of them).

\section{Algebra of gain graphs and their states} \label{states}

In applications of the kind in Section \ref{applied} gain graphs are normally used in the following setup. Let
$(\Gamma ,g,G)$ be a gain graph.  Let $Q$ be a set, which we call the set of \emph{qualities}, and suppose $G$
acts on $Q$; for $q\in Q$ and $h\in G$ we denote the result of the action by $[q]h$.  For instance, it may be
that $Q=G$ and the action is by right multiplication (the \emph{right regular action}). The combination
$(\G,g,G,Q)$ is a \emph{permutation gain graph}.  (The idea was introduced by Gross \& Tucker (1977) under the
name ``permutation voltage graph''.)  A \emph{state} of the graph is an assignment $s$ of qualities to the
vertices of $\Gamma$. We want to find a state such that for any two adjacent vertices $v$ and $w$ their
qualities $s(v)$ and $s(w)$ are connected by $s(w)=[s(v)]g(vw)$. Such a state is called \emph{satisfied} with
respect to $g$, or \emph{$g$-satisfied}. The set of all $g$-satisfied states is denoted by $\Sat_g(\G,G,Q)$.
The set of all states that are satisfied with respect to at least one element of $G$ is denoted by
$\Sat(\G,G,Q)$. If $H$ is a subgroup of $G^E$, then $\Sat_H(\G,G,Q)$ is the set of all states that are
satisfied with respect to some $g \in H$. A $g$-satisfied state always exists if $g$ is balanced, because in
this case one can arbitrarily choose the state of one vertex in each connected component, then assign the
states of other vertices according to the rule of satisfaction. The exact rule is that, if $s(v_0)$ is the
arbitrary state and $v$ is a vertex connected to $v_0$ by a walk $W$, then $s(v) = [s(v_0)]g(W)$.

The set $G^E$ of all gain maps with values in $G$ forms a group under
componentwise multiplication, $gg'(e)=g(e)g'(e)$.
If $G$ is abelian the set $\Bal(\G,G)$ of all balanced gain maps is a
subgroup of $G^E$.
If, in addition, $\G$ is connected, then $|\Bal(\G,G)| = |G|^{|V|-1}$,
and $\Bal(\G,G) \cong G^{|V|-1}$ if $G$ is abelian, because one can
assign gains arbitrarily to the edges of a spanning tree and then other
edge gains are uniquely determined by balance.
In general, $\Bal(\G,G) \not\cong G^{|V|-1}$ if $G$ is nonabelian.
From now on \emph{we assume that $\G$ is connected.}

If $Q$ is a group, the set $Q^V$ of states forms a group under
componentwise multiplication, $ss'(v)=s(v)s'(v)$.
From now on \emph{we assume that $Q$ is a group.}  We do not assume, in
general, that the action of $G$ on $Q$ is a group action.

If $G$ does act on $Q$ in the manner of a \emph{group action} (i.e.,
through a homomorphism $G \to \Aut Q$), then the set of all $g$-satisfied
states, $\Sat_g(\G,G,Q)$, forms a subgroup of $Q^V$.
If in addition $\G$ is connected and $g$ is balanced, then
$$
\Sat_g(\G,G,Q) \cong Q
$$
because the choice of $s(v)$ at one vertex $v$ determines $s$.

Another way $G$ may act on $Q$ is by a \emph{product action}, that is, where
\begin{equation}
[qq'](ff')=([q]f)([q']f'). \label{E:productaction}
\end{equation}
Affine action of $\R^n$ on itself (by translations) is an example of a
product action in  geometry.
A product action is very restricted.

\begin{lemma} \label{L:productaction}
An action of $G$ on $Q$ is a product action if and only if it has the
form $[q]f=\sigma(q)c(f)$, where $\sigma: Q \to Q$ is an endomorphism of
$Q$ and $c: G \rightarrow Q$ is a homomorphism of $G$ into the center of
$Q$.
\end{lemma}

\begin{proof}
\comment{\emph{NEW COMMENT}:  Can this proof be omitted?  ONLY IN CASE
WE ARE ASKED FOR FURTHER REDUCTIONS.}
Any action of this form is clearly a product action. For the converse, compare the two
sides of \eqref{E:productaction}. First, with $q=q'=1$; this establishes the existence of a homomorphism $c$
given by $c(f)=[1]f$. Next, with $q'=1$, and $f=1$; this establishes that $[q]f'=([q]1)\cdot c(f')$. Taking
$q=1$ and $f'=1$ establishes that $[q]f=c(f)q'$; combining this with the previous conclusion we get $[q]1\cdot
c(f)=[q]f=c(f)\cdot[q]1$ for all $q \in Q$. So, $c(f)$ is in the center of $Q$. Taking $f=f'=1$ shows that
$\sigma: q \mapsto [q]1$ is a homomorphism $Q \to Q$.
\end{proof}

We call  $(g,s)\in G^E \times Q^V$ a \emph{satisfied pair} if $s$ is a
satisfied state with respect to the gain map $g$.
With a product action of $G$ on $Q$, the set of satisfied pairs forms a
group $\SatPair(\G,G,Q)$ with respect to  componentwise multiplication,
$(g,s)(g',s')=(gg',ss')$.
If $H \subseteq G^E$, then $\SatPair_H(\G,G,Q)$ stands for the set of all
satisfied pairs $(s,g)$ where $g \in H$.
If $G$ is abelian and $H$ happens to be a subgroup of $\Bal(\G,G)$,
$\SatPair_H(\G,G,Q)$ forms a subgroup of $\SatPair(\G,G,Q)$.

We call an action of $G$ on $Q$ \emph{fixed-point free} if every element of $G$, other than the identity, acts
without fixed points. If $\G$ is connected and the action of $G$ is fixed-point free, then $\Sat(\G,G,Q) \cong
\Bal(\G,G) \times Q$,  because, as shown in the following lemma, a gain mapping $g \in G^E$ and the value $q
\in Q$ at one vertex together determine a satisfied state if and only if $g \in \Bal(\G,G)$.
Ffrom the definitions we have:

\begin{lemma}\label{L:FixedPointFreeAction}
Suppose the action of $G$ on $Q$ is fixed-point free.  Then a permutation
gain graph $(\G,g,G,Q)$ has a satisfied state if and only if $(\G,g,G)$
is balanced, and then the state is determined by its value at any one vertex.
\end{lemma}
\begin{proof}
Given $(\G,g,G,Q)$, let $s$ be a state.  Let $W_v$ denote a simple closed walk from $v$
to $v$.  Then $s$ is satisfied if and only if $[s(v)]g(W_v) = s(v)$ for all choices of
$v$ and $W_v$.  Under the hypothesis on fixed points, this is true if and only if
$g(W_v) = 1$ for all choices of $W_v$, i.e., if and only if $g$ is balanced.
\end{proof}

\comment{This remark can be removed if we are asked for further reductions.}
On the other hand, if the action if not fixed-point free and $\G$ is
not a tree, then one can always define an unbalanced gain map that has
a satisfied state. The lemma implies, in particular, that if $H$ is a
subgroup of $\Bal(\G,G)$, then there is a natural bijection between
$\Sat_H(\G,G,Q)$ and $H \times Q$.

Suppose $G$ acts on $Q$ in the manner of a group action. If $Q$ is a
module over a ring $\cR$ and the group action is linear over $\cR$,
then for any gain map $g$,  $\Sat_g(\G,G,Q)$ is also an $\cR$-module.
Furthermore, by Lemma \ref{L:FixedPointFreeAction}
\begin{equation} \label{E:Sat_g=Q}
\Sat_g(\G,G,Q) \cong Q.
\end{equation}

If $G$ is also an $\cR$-module, then the group of gain mappings, $G^E$,
is an $\cR$-module.
If, moreover, $H$ is a submodule of $\Bal(\G,G)$, then the action of $G$
on $Q$ is indeed a product action, $\SatPair_H(\G,G,Q)$ is a subgroup of
$\SatPair(\G,G,Q)$, and
\begin{equation} \label{E:SatPair_H=Q+H}
\SatPair_H(\G,G,Q) \cong  Q \oplus H.
\end{equation}
In particular, $\SatPair(\G,G,Q) \cong  Q \oplus G^{|V|-1}$.

If the action of module $G$ on module $Q$ is also fixed-point free, then
by Lemma \ref{L:FixedPointFreeAction},
\begin{equation}\label{E:Sat_H=Q+H}
\Sat_H(\G,G,Q) \cong \SatPair_H(\G,G,Q) \cong  Q \oplus H .
\end{equation}

A special case of permutation gain graphs has been studied by Joswig (2002) and
Izmestiev \& Joswig (2003) in connection with coloring of tilings of manifolds and
problems in the geometric topology of covering spaces.  Let $\G$ be the dual graph of a
finite, pure, $d$-dimensional, strongly connected simplicial complex $\Delta$ on $n$
vertices, and let $G$, the gain group, be $S_n$, the symmetric group on the vertices of
$\Delta$.  The gain mapping $g$ takes an edge connecting adjacent $d$-simplices whose
symmetric difference is a pair $\{i,j\}$ to the transposition that exchanges $i$ and
$j$ (Joswig calls this a ``perspectivity'').  Joswig calls the image $g(\pi(\G))$ of
the group of closed edge-walks beginning and ending at a vertex $v \in \G$ the group of
``projectivities'' and writes it $\Pi(\Delta,v)$.  The main object of Izmestiev and
Joswig's studies is the size and structure of this group for various combinatorial
manifolds. See Appendix for more details on projectivities and their applications in
combinatorial topology.

\section{Gain graphs in geometry} \label{applied}

Given a tiling of $\R^d$ by convex polyhedra, it is natural to ask
whether this tiling can be thought of as the vertical projection of a
convex PL-surface or, more generally, any PL-surface.  A PL-function on
$\R^d$ is known as a \emph{$C_1^{0}$-spline} in approximation theory.
The problem of finding the dimension of the space of such splines on a
given tiling has received significant attention from discrete geometers
and spline enthusiasts (see, e.g., Whiteley (1996) or Billera (1989) for
details).
Voronoi (1908) showed that any \emph{simple} tiling $T$, i.e., such that
each vertex is incident to exactly $d+1$ tiles, is indeed the projection
of a convex surface, which he called a \emph{generatrix} of $T$;
furthermore, a surface that gives this tiling as a projection is uniquely
determined by the choice of an affine support function and a dihedral
angle.

A tiling is an example of a more general concept, namely a PL-realization
$M$ in $\R^{d}$ of a $d$-dimensional manifold $\cM$.
In this section, as an application of the Binary Cycle Test, we ask about
possible PL-realizations of $\cM$ in $\R^{d+1}$ that give $M$ as a
vertical projection.

If $M$ is a PL-realization in $\R^{d}$ of a cell-complex $\cM$
(\emph{cell-complex} defined, e.g., in the sense of Seifert \& Threlfall
(1980); the PL-realization is what is called a \emph{geometric
cell-complex}),
and $L$ is another PL-realization of $\cM$, but in $\R^{d+1}$,
such that for each cell $C$ of $\cM$ its realization in $M$ is the
vertical projection of its realization in $L$, then $L$ is called a
\emph{lifting} of $M$.
The notion of lifting naturally generalizes that of $C_1^{0}$-spline.
The graph of a $C_1^{0}$-spline is a lifting of the tiling over which it
is defined, while a lifting of a PL-realization $M$ in $\R^{d}$ of a
$d$-manifold $\cM$ is the graph of a $C_1^{0}$-spline if and only if $M$
is an embedding.
Obviously, liftings can be added, just like $C_1^{0}$-splines. The space
of liftings of $M$ is denoted by $\Lift(M)$.
A trivial kind of lifting is to a hyperplane in $\R^{d+1}$.  The opposite
kind is a sharp lifting: we call a PL-realization in $\R^{d+1}$ of a
$d$-complex \emph{sharp} if any two adjacent $d$-cells are realized on
different affine hyperplanes.

Of course, questions about projections of piecewise-linear surfaces
belong to the realm of projective geometry.  To make our presentation
more visual and accessible we do all our geometry in the real affine
space $\R^d$, although at the end of the section we give our result a
proper projective interpretation.

Although some of the constructions that we treat below were considered in Rybnikov
(1999), the Binary Cycle Test is formulated incorrectly there (in Lemma 4.1).  Roughly
speaking, that formulation does not have any restriction on gain groups, due to the
fact that the author had in mind only abelian groups. However, since the test was
applied there for abelian, torsion-free
groups only, all the geometric results of Rybnikov (1999) 
about PL-realizations of (and splines over) a manifold $\cM$ with
$H_1(\cM;\Z_2)=0$ hold true if one assumes that the manifold has finitely
many cells; this last requirement was also omitted in Rybnikov (1999).
(When $\cM$ is simply connected the finiteness restriction is not
necessary, since in this case there is an alternative argument, which
uses the fundamental group---see Ryshkov \& Rybnikov (1997) for proofs
and examples.)

\subsection{Technical Preparation}

We use standard notions of PL-topology.  (See, e.g., Munkres (1984) or Seifert \&
Threlfall (1980).) A \emph{regular CW-complex} $\cK$ is a locally finite CW-complex in
which all gluing mappings are homeomorphisms. A regular CW-complex has a barycentric
subdivision $\cK^{\flat}$, which is a simplicial complex. $\cK$ also has a natural dual
structure, called the \emph{dual cell-decomposition}; we will denote the dual
CW-complex by $\cK^*$.  Note that $\cK$, $\cK^*$ and $\cK^{\flat}$ have the same
support and are homeomorphic. The face lattices of $\cK$ and $\cK^*$ are dual. A cell
is always assumed to be \emph{relatively open }and the closure of a cell need not be
homeomorphic to a ball. A $k$-cell $C^*$ of $\cK^*$ is the union of all simplices of
$\cK^{\flat}$ whose closures share the barycenter of a $(d-k)$-cell $C$ of $\cK$. The
boundary of $C^*$ is called the link of $C$. The dimension of $\cK$ is the maximum
dimension of any cell of $\cK$. \emph{Facets} of $\cK$ are cells of dimension
$\dim\cK-1$, and \emph{ridges} of $\cK$ are cells of dimension $\dim\cK-2$. Two
$k$-cells are \emph{adjacent} if their closures share a $(k-1)$-cell. A pure
$d$-dimensional CW-complex is \emph{strongly connected} if any two $d$-cells can be
joined by a sequence of adjacent $d$-cells. We use $\textrm{C}^k (\cK,A)$ and
$\textrm{C}_k (\cK,A)$ to refer to the $\R$-spaces of $k$-co-chains and $k$-chains over
a real vector space $A$.

A \emph{homology sphere} is defined inductively: a homology 0-sphere is a
standard 0-sphere, and a higher-dimensional homology sphere is a strongly
connected regular CW-complex such that: (1) all its homology groups are
like those of a standard sphere and (2) the boundaries of all cells of the
dual cell-decomposition are homology spheres.  A homology ball is a cone
over a homology sphere. A finite \emph{homology $d$-manifold (with
boundary)} is a strongly connected regular CW-complex with finitely many
cells, in which the link of each $k$-cell, for $0 \leq k < d$, is a
homology sphere (or ball, if the cell belongs to the boundary) of
dimension $d-k-1$.

Seifert \& Threlfall (1980) allow the cells of a homology manifold to be
homology balls; then the regular CW-complexes in our definitions should be
replaced by regular homology CW-complexes. All of our results hold for
this more general interpretation.

If $\cK$ is a regular CW-complex and $r$ is a continuous map from the
barycentric subdivision  $\cK^{\flat}$ to $\R^N$ such that (1) each
abstract simplex is mapped  to an affine simplex, and (2) the simplices
of $\cK^{\flat}$ making up a $k$-cell of $\cK$ are mapped to the
same affine $k$-flat of $\R^N$, then $K=(\cK,r)$ is called a
\emph{PL-realization} of $\cK$.  The realization is called
\emph{non-degenerate} if $r$, restricted to the closure of each cell of
$\cK$, is a homeomorphism.

\subsection{Local Considerations}

An \emph{oriented facet} is an ordered pair of adjacent $d$-cells.
The minimum possible number of $d$-cells that make contact at a $k$-cell
of a convex tiling of $\R^{d}$, or of a $(d+1)$-polytope, is $d-k+1$.
We call the star of a $k$-cell in a manifold \emph{simple} (a.k.a.\
\emph{primitive} in the theory of parallelohedra) if the cell is incident
to only $d-k+1$ $d$-cells.
A $d$-manifold is called \emph{$k$-simple} if the star of each $k$-cell
is simple.  A manifold with boundary is $k$-simple if the stars of all
its interior $k$-cells are simple.  Notice that $k$-simplicity implies
$m$-simplicity for $m>k$.

The \emph{dual graph} $\G(\cM)$ of a homology $d$-manifold $\cM$
(possibly with boundary) is the graph of adjacency of $d$-cells.  It is
the 1-skeleton of the dual cell-decomposition $\cM^*$.
The vertices are the $d$-cells of $\cM$ and the edges are the pairs of
$d$-cells that share a common facet.
The edges of $\G(\cM)$ can be thought of as the internal facets of $\cM$.
In the dual graph of a $k$-simple manifold, the subgraphs corresponding
to the stars of $k$-cells are complete subgraphs on $d-k+1$ vertices.
Subgraphs of $\G(\cM)$ that correspond to the stars of ridges are
circles; we call them \emph{ridge circles} and, in the case when the
stars of ridges are simple, \emph{ridge triangles}.

A \emph{reciprocal} for a PL-realization $M$ in $\R^{d}$ of a homology
manifold $\cM$ is a rectilinear realization $R$ in $\R^{d}$ of the dual
graph $\G(\cM)$ such that the edges of $R$ are perpendicular to the
corresponding facets. If none of the edges of $R$ collapse into a point, the
reciprocal is called \emph{non-degenerate}; a one-point reciprocal is
called \emph{trivial}.  A reciprocal for the star of a cell $C$ is called
a \emph{local reciprocal} of $C$. Reciprocals were originally considered
by Maxwell (1864) and Cremona (1890) in connection with stresses in plane
frameworks. Crapo \& Whiteley (1993) gave a modern treatment of the
Maxwell--Cremona theory of reciprocals, stresses, and liftings for
2-manifolds. Obviously, reciprocals can be added and multiplied by real
scalars.  Geometrically, multiplication corresponds to scaling. The
addition is induced by vector addition of the vertices. Thus, for a given
$M$ the reciprocals form a linear space. This space, factored by all
translations of $\R^{d}$, is denoted $\Rec(M)$.

A PL-realization of a homology $d$-manifold is called \emph{generic} if
no pair of adjacent facets lies in the same $(d-1)$-plane.

\begin{lemma}\label{Uniqueness of Local Reciprocal}
For each non-degenerate  generic PL-realization in $\R^{d}$ of a simple
star of a cell there is a
non-degenerate reciprocal.  It is unique up to translations, homotheties,
and central inversions.
\end{lemma}

\begin{proof}
The existence of a non-degenerate reciprocal for a simple star is a
trivial exercise in linear algebra.  Such a reciprocal is the 1-skeleton
of a simplex whose edges are perpendicular to the corresponding facets.
Since directions of edges and the angles between them are uniquely determined
by the geometry of the star, the reciprocal is unique up to affine
transformations preserving directions and angles.
\end{proof}

\begin{lemma}\label{lifting simple star}
If $Star$ is a non-degenerate generic PL-realization of a simple star of
a cell, then $\dim\Lift(Star)=d+2$.
\end{lemma}

\begin{proof}
Any three $d$-cells of $Star$ form a triangle in the dual graph.
Choose two $d$-cells in $Star$.
By fixing two affine functions that lift the union of
these $d$-cells  to $\R^{d+1}$, so that these functions coincide on the
common facet of the cells,  we fix the lifting for all of $Star$.  Since
a lifting of two adjacent $d$-cells is defined by $d+2$ parameters (e.g.,
the coefficients of an affine function in $d$ variables and a dihedral
angle in $(d+1)$-space), $\dim\Lift(Star)=d+2$.
\end{proof}

\subsection{From Local to Global}

In the remainder of the section we show that, subject to the topological
restriction that $H_1(\cM;\Z_2)=0$ and the combinatorial restriction of
$(d-3)$-simplicity, liftings of $(d-2)$-simple manifolds behave in the
same way as those of simple stars.  (Although a tiling of $\R^d$ may not
be a CW-complex, all our theorems hold for the case when $M$ is a finite
tiling.)

The following lemma is a homological version of what some call
Poincar\'e's $(d-2)$-face principle (1882) in discrete geometry, which,
basically, says that for a group of isometries to tile a
Euclidean space, sphere, or hyperbolic space with copies of a
$d$-polyhedron $T$ it is enough that for each $(d-2)$-face $F$ of $T$
there are group elements $g_1,\ldots,g_k$ such that
$T,g_1(T),\ldots,g_k(T)$ fill the space, face-to-face, around $F$ without
gaps or overlaps.  See Dolbilin (2000) for an exact formulation and recent
generalization of this principle.

\begin{lemma}
\label{Balance Condition For Manifolds}
Let $\cM$ be a finite homology $d$-manifold with $d \geq 2$ and with
$H_1(\cM;\Z_2)=0$.  Let $(\G(\cM),g,G)$ be an abelian gain graph whose
underlying graph is $\G(\cM)$ and whose gain group has no odd torsion.
Then the gain graph is balanced if and only if the ridge circles of
$\G(\cM)$ are balanced.
\end{lemma}

\begin{proof}
The ``only if'' part follows from the definition of balanced graph.

Consider the dual cell-decomposition $\cM^*$ of $\cM$.  The 2-cells of
this decomposition are in one-to-one correspondence with the ridges of
$\cM$.  Since $H_1(\cM;\Z_2)=0$ and the underlying spaces of $\cM$ and
$\cM^*$ are the same, the boundaries of the 2-cells of $\cM^*$ generate
$Z_1(\cM^*;\Z_2)$. These boundaries are circles and they are balanced by
our assumption. According to the Binary Cycle Test for finite graphs,
$(\G(\cM),g,G)$ must be balanced.
\end{proof}

For instance, Lemma \ref{Balance Condition For Manifolds} applies to
homology spheres of dimension at least 2.

\begin{theorem}
\label{global reciprocal}
Let $\cM$ be a finite $(d-3)$-simple homology $d$-manifold without
boundary, with $d \geq 3$, and with $H_1(\cM;\Z_2)=0$.
Let $M$ be a non-degenerate generic PL-realization of $\cM$ in $\R^{d}$.
Then $M$ has a non-degenerate reciprocal and $\dim \Rec(M)=1$.
\end{theorem}

\begin{proof}
The main ingredients of the proof are the dual cell-decomposition $\cM^*$
of $\cM$, the notion of gain graph, and the Binary Cycle Test for balance.

Let $\G=\G(V,E)$ be the dual graph of $\cM$.  An oriented facet can be
thought of as an oriented edge of $\G$.

From the algebraic standpoint a reciprocal for $M$ is a linear operator $r$ from the
vector space $(\bbR^d)^V$ of 0-cochains on $\G(V,E)$ with values in $\R^d$, which has
the property that $\delta r(C_i,C_j)$ is orthogonal to the common facet $F_{ij}$ of
$C_i$ and $C_j$, where $\delta: C^0(\G,\R^d) \to C^1(\G,\R^d)$ is the coboundary map,
defined by $\delta r(C_i,C_j) = r(C_j)-r(C_i)$ for each oriented facet $(C_i,C_j)$.
(See Appendix for more detailed description of the above.) Thus, the edges
$\{(r(C_i),r(C_j))\}$ of a non-degenerate reciprocal have a natural interpretation as
the linear parts of the equations defining the facets of $M$.  These linear parts are
often called the \emph{covectors} of the facets of $M$.  It follows from the definition
that
\begin{equation} \label{vector gains}
\delta r(C_1,C_2)=-\delta r(C_2,C_1)
\end{equation}
and that
\begin{equation} \label{vector ridge balance}
\delta r(C_1,C_2)+\delta r(C_2,C_3)+\delta r(C_3,C_1) = 0
\end{equation}
whenever $C_1$, $C_2$, $C_3$ share a common ridge.  Equation \ref{vector
gains} means that $\delta r$ is a gain map on $\G$ with gain group $\R^d$.
The gain graph is balanced because $\delta r$ is a coboundary.
Because $\G$ is connected, $r$ is completely determined by $\delta r$ and
one point $r(C)$; indeed, $r$ is any satisfied state of
$(\G,\delta r,\R^d,\R^d)$, where $\R^d$ acts on $\R^d$ by translation
(that is, the right regular action).  Thus,
\begin{equation}
\label{reciprocals as coboundaries}
\Rec(M) = \{\delta r : r \text{ is a reciprocal for } M\},
\end{equation}
the space of coboundaries of reciprocals.

We associate with each oriented facet $(C_1,C_2)$ a unit vector
$\n(C_1,C_2)$, normal to the common facet of $C_1$ and $C_2$, so that
$\n(C_1,C_2)=-\n(C_2,C_1)$.  Henceforth, these normals are fixed.

Consider the \emph{facet graph} $\Phi$, whose vertices are the (unoriented)
facets of $\cM$ and whose edges are the pairs of facets sharing ridges.
Orient the edges of $\Phi$ in an arbitrary way.
Take a ridge $R$ where $d$-cells $C_1$, $C_2$, $C_3$ make contact.
Denote by $F_{ij}$  the common facet of $C_{i}$ and $C_{j}$.
Up to scaling there is only one nontrivial linear combination of
$\n(C_1,C_2)$, $\n(C_2,C_3)$, $\n(C_3,C_1)$ equal to zero.  Let
\begin{equation*}
\alpha_{12}\n(C_1,C_2) + \alpha_{23}\n(C_2,C_3) + \alpha_{31}\n(C_3,C_1) = 0
\end{equation*}
be one such combination with nonzero coefficients; it exists because $M$
is a non-degenerate PL-realization. For the oriented edge $F_{12}F_{23}$
in the facet graph, define
\begin{equation*}
g(F_{12}F_{23}) = \alpha_{23}/\alpha_{12} .
\end{equation*}
In this way we define the map $g$ on all oriented edges of the facet
graph.  We call $(\Phi,g,\R^*)$ the \emph{facet gain graph}.

\begin{lemma}
\label{balance of facet gain graph}
$(\Phi,g,\R^*)$ is a balanced gain graph.
\end{lemma}

\begin{proof}
Obviously, $g(F_{12}F_{23})$ depends only on the geometric realization
$M$ of $\cM$ and on the (fixed) choice of normals $\n(A,B)$ for oriented
facets.  It is clear that $g(F_{12}F_{23})=g(F_{23}F_{12})^{-1}$ and
\begin{equation}
\label{product over cycle of facet graph}
g(F_{12}F_{23}) \cdot g(F_{23}F_{31}) \cdot g(F_{31}F_{12})=1 ;
\end{equation}
that is, $g$ is a gain map on $\Phi$ and each triangle corresponding to a
ridge (a \emph{ridge triangle of $\Phi$}) is balanced.  Furthermore,
\begin{equation}
\label{key equation}
\n(C_1,C_2) + g(F_{12}F_{23})\n(C_2,C_3) +
g(F_{23}F_{31})g(F_{12}F_{23})\n(C_3,C_1) = 0.
\end{equation}

The dual cell-decomposition $\cM^*$ of $\cM$ can be refined in a useful way.
Each 2-cell of $\cM^*$ is a triangle.
Let $\cM^{\dag}$ be the cell-complex obtained from $\cM^*$ by subdividing
each edge in the middle and then joining the \emph{new vertices} by three
\emph{new edges} within each 2-cell, thereby partitioning the 2-cell into
four triangles.

Since $\cM^*$ and $\cM^\dag$ are identical to $\cM$ as topological spaces,
$$
H_1(\cM^\dag;\Z_2)=H_1(\cM^*;\Z_2)=H_1(\cM;\Z_2)=0.
$$
Moreover, $\Phi$ is a subgraph of $\Sk_1(\cM^\dag)$.  Consequently, for
any $c\in Z_1(\Phi;\Z_2)$ there is a 2-chain
$\Delta \in C_2(\cM^\dag;\Z_2)$ such that $c=\partial \Delta$. The chain
$\Delta$ can be decomposed as $\Delta_1+\Delta_2$, the sum of
$\Delta_1 = \sum \Delta_{1j}$, where each $\Delta_{1j}$ is a triangle
whose vertices are new vertices corresponding to facets of $\cM$ making
contact at a ridge, and $\Delta_2=\sum \Delta_{2i}$, where each
$\Delta_{2i}$ has one vertex in $\Sk_0(\cM^*)$, corresponding to a cell
of $\cM$, and two vertices that are new vertices corresponding to facets
of that cell. For each $d$-cell $C$ of $\cM$, let $\Delta_{2C}$ be the
sum of all the $\Delta_{2i}$ having as one vertex that vertex of $\cM^*$
that corresponds to $C$. Then $c=\sum \partial \Delta_{1j}+\sum \partial
\Delta_{2C}$ where each term in the first sum is a ridge triangle of
$\Phi$ and each term in the second sum is a cycle on the boundary of a
$d$-cell of $\cM$.

This shows that ridge triangles and cycles lying on the boundaries of
$d$-cells generate $Z_1(\Phi;\Z_2)$. We prove next that $(\Phi,g,\R^*)$
is balanced. Since $\cM$ is $(d-3)$-simple, the boundary of each $d$-cell
$C$ is $((d-1)-2)$-simple.  $\Phi_C$, the subgraph of $\Phi$ on the
vertices that belong to $C$, is the dual graph $\G(\partial C)$ of the
boundary cell-complex of $C$. The corresponding gain subgraph is
$(\Phi_{C},g,\R^*)$. By Equation \eqref{product over cycle of facet
graph} the ridge circles of $\Phi_{C}$ are balanced. Since
$\dim\partial C\geq 2$, Lemma \ref{Balance Condition For Manifolds}
applied to $\partial C$ implies that $(\Phi_{C},g,\R^*)$ is balanced and,
therefore, all cyclic orientations of binary cycles of $\Phi_{C}$ have
gain 1. We have seen that binary cycles of these two types generate
$Z_1(\Phi;\Z_2)$. By the Binary Cycle Test for finite graphs,
$(\Phi,g,\R^*)$ is balanced.
\end{proof}

Suppose the set of qualities $Q$ is  $\R$ and $\R^*$ acts on $\R$ by
multiplication. The action of $\R^*$ on $\R\setminus0$ is fixed-point
free and therefore, by Lemma \ref{L:FixedPointFreeAction} and the obvious
properties of the all-zero state, $\Sat_g(\Phi,\R^*,\R) \cong \R$. We are
now going to explain why $\Sat_g(\G,\R^*,\R)$ can be regarded as a
1-dimensional subspace of $(\R^d)^E$, where $E$ is the set of edges of $\G$.

Fix a non-zero satisfied state $s$ of $(\Phi,g,\R^*,\R)$.  We denote by
$s(C_i,C_j)$ the value of $s$ on a vertex of $\Phi$ corresponding to a
facet between $d$-cells $C_i$ and $C_j$.
We will construct yet another gain graph, $(\G,h_s,\R^{d})$, this time
with the gain group $\R^{d}$.
Assign to each oriented edge $C_1C_2$ of $\G$ the vector $h_s(C_1C_2) =
s(C_1,C_2)\n(C_1,C_2)$.

\begin{lemma}\label{balance of dual gain graph}
$(\G,h_s,\R^{d})$ is a balanced gain graph.
\end{lemma}

\begin{proof}
Since $\n(C_1,C_2)=-\n(C_1,C_2)$ and $s(C_1,C_2)=s(C_2,C_1)$,
$h_s(C_2C_1)=-h(C_1C_2)$. Since $H_1(\cM^*;\Z_2) = H_1(\cM;\Z_2) = 0$, by
Lemma \ref{Balance Condition For Manifolds} balance of the ridge triangles
of $\G$ implies balance of the gain graph.  To conclude that these
triangles are balanced we have to show that
$$
s(C_1,C_2)\n(C_1,C_2)+s(C_2,C_3)\n(C_2,C_3)+s(C_3,C_1)\n(C_3,C_1) = 0.
$$
Let us divide this equation by $s(C_1,C_2)$ and recall that (with indices
taken modulo 3)
$s(C_{i+1},C_{i+2})/s(C_{i},C_{i+1}) = g(F_{i,i+1}F_{i+1,i+2})$ for the
star of each ridge.  So, we need
\begin{equation*}
\n(C_1,C_2)+g(F_{12}F_{23})\n(C_2,C_3)+g(F_{12}F_{31})\n(C_3,C_1) = 0.
\end{equation*}
This is equivalent to Equation \eqref{key equation}. Thus each ridge
triangle is balanced, and $(\G,h_s,\R^{d})$ is balanced by the Binary Cycle
Test.
\end{proof}

We want to prove that the elements of $\Sat_g(\Phi,\R^*,\R)$ correspond
to the reciprocals (modulo $\R^d$) in one-to-one fashion; hence
$\dim\Rec M = 1$.
Our strategy is to show that $\Sat_g(\Phi,\R^*,\R) \cong \Rec(M)$ as real
vector spaces by producing a specific isomorphism $\rho$. We define
$\rho(s) = h_s$.  It is obvious that $\rho$ is an injection and a
vector-space homomorphism. Let $H$ be the image of $\rho$. Then
$H \subseteq \Bal(\G,\R^d)$ and $\dim H = 1$.

Consider now the permutation gain graph $(\G,h_s,\R^{d},\R^{d})$, where
the first three components are as above and $\R^{d}$ acts on $\R^{d}$ by
translation.  This action is obviously fixed-point free. Since a
reciprocal is completely defined by the coordinates of its vertices, all
reciprocals form a linear subspace of $(\R^d)^{V}$.  A satisfied state for
$(\G,h_s,\R^{d},\R^{d})$ is a reciprocal $r$ whose edge lengths
have been decided by the gain graph $(\G,h_s,\R^{d})$.  That is, since
$(\G,h_s,\R^{d})$ is balanced and connected, we can start from any point
in $\R^d$ as a vertex of the reciprocal and construct all other vertices.
Then $h=\delta r$.
By Lemma \ref{L:FixedPointFreeAction} this construction gives all of
$\Sat_{h_s}(\G,\R^{d},\R^{d})$; varying $s$, it gives all of
$\Sat_H(\G,\R^{d},\R^{d})$.
Furthermore, it shows that $H \subseteq \{\delta r \: | \: r \textrm{ is a
reciprocal} \}$, the reciprocal coboundary
space, which is $\Rec M$ by Equation \ref{reciprocals as coboundaries}.
The non-zero elements of $H$ \now{correspond (bijectively) to} the
non-degenerate reciprocals.

Conversely, the edges of a reciprocal $r$ determine a balanced gain map
$\delta r$ with values in the gain group $\R^d$.  By the orthogonality
property of a reciprocal and Equation \eqref{vector ridge balance},
$\delta r$ has the form $\rho(s)$ for some non-zero
$s \in (\Phi,g,\R^*,\R)$. Therefore,
$H = \{\delta r \: | \: r \textrm{ is a reciprocal} \} =\Rec(M)$.

We established earlier that $\dim_{\R} H=1$.
It follows that $\dim_{\R} \Rec(M)=1$.
\end{proof}

\begin{theorem}
\label{main lift theorem}
Let $\cM$ and $M$ be as in Theorem \ref{global reciprocal}.  Then $M$ has
a sharp lifting and $\dim_{\R} \Lift(M)$ is $d+2$.
\end{theorem}

\begin{proof}
Since $\cM$ is $(d-2)$-simple, in any lifting of $M$ the position of a
$d$-cell $C$  from the star $\Star (R)$ of a ridge $R$ is completely
defined by the positions of the other two $d$-cells of $\Star (R)$.
Therefore, if the lifting of $\Star (R)$ can be extended to  a lifting of
$M$, the resulting lifting of $M$ is completely defined by the lifting of
$\Star (R)$.  By Lemma \ref{lifting simple star},
$\dim \Lift (M) \le d+2$.

Let us now prove that $\dim \Lift (M) = d+2$. Let $H$ be as in the above
theorem, i.e., let $H$ be the subspace of $\Bal(\G,\R^d,\R^d)$ that
consists of those realizations of the dual graph $\G$ of $\cM$ that are
reciprocals of $M$. Consider gain graphs where the underlying graph is
the dual graph $\G$ and the gain group is the additive group $A_d$ of all
affine functions on $\R^d$ (this is isomorphic to $\R^{d+1}$). $A_d$
acts on itself by addition and this action is fixed-point free.
Therefore, by Lemma \ref{L:FixedPointFreeAction}, the liftings of $M$ can be
identified with all satisfied states of permutation gain graphs
$(\G,h_1,\R^{d+1},A_d)$ where $h_1$ is taken from the subgroup $H_1$ of
$\Bal(\G,\R^{d+1})$ that consists of gain maps of the form
$h_1(F)=\e_F \cdot (\x-\c_F)$, where $F$ is an oriented facet, $\e_F$ is the
corresponding oriented edge in some reciprocal, $\c_F$ is any point on
$F$, and $\x$ is the variable vector. This gain map can be written as
$h_1(F)=h(F)+c(F)$, where $h(F)=\e_F \cdot \x$ is a gain map in
$H \subseteq \Bal(\G,\R^{d})$. We know that $H \cong \R$. Therefore, a
state $s$ of $(\G,h_1,\R^{d+1},A_d)$ corresponds to a lifting if and only
if the gain map $c(F)=-\e_F \cdot \x $ is balanced. It is easy to see
that all cycles of $\G$ corresponding to ridges of $\cM$ are balanced with
respect to this gain map. These cycles generate all of
$H_1(\cM^*;\Z_2)=H_1(\cM;\Z_2)$. Therefore, by the Binary Cycle Test $c$
is balanced. Thus, by Lemma \ref{E:Sat_H=Q+H},
$\Lift(M) \cong A_d \oplus \R \cong \R^{d+2}$.
\end{proof}

\begin{corollary}
Let $\cM$ be a finite simple homology $d$-manifold, where $d \geq 3$,
such that $H_1(\cM,\Z_2)=0$. Let $M$ be a non-degenerate generic
PL-realization in $\R^{d}$ whose $d$-cells are convex.  Then $M$ is a
projection of a convex polyhedron in $\R^{d+1}$.
\end{corollary}

\begin{proof}
Since $\cM$ is $(d-3)$-simple, there is a sharp lifting $L$ of $M$ to
$\R^{d+1}$.  Since all cells of $\cM$ are homeomorphic to convex
polytopes, and there are $d-k+1$ $d$-cells making contact at each
$k$-face, $\cM$ is a topological manifold. The sharp lifting $L$ is then a
simple PL-surface in $\R^{d+1}$, realizing abstract manifold $\cM$.  $L$
must be an immersion, since each cell is embedded and the star of each
vertex is embedded. So, we have a realization of a $d$-manifold in
$\R^{d+1}$ which has the following properties: (1) it is an immersion, (2) it
is locally (non-strictly) convex at each point, and (3) it is strictly
convex at each vertex. By Van Heijenoort's (1952) theorem $L$ is the
boundary of a convex body. Since  $\cM$ is finite and $L$ is
piecewise-linear, $L$ is a convex polyhedron.
\end{proof}

Although Theorem \ref{main lift theorem} is formulated for vertical
projections, it, indeed, holds for a central projection from $\R^{d+1}$ to
$\R^{d} \subset \R^{d+1}$.  It also holds for a central projection from
the origin of $d$-manifolds realized in $\R^{d+1}$ to a $d$-sphere in
$\R^{d+1}$ (see Rybnikov, 1999), or projection from $(0,...,0,-1)$ of
$d$-manifolds realized in $\R^{d+1}$ to the Minkowski hyperboloid in
$\R^{d+1}$. Of course, the reason for such universality is that the
theorem is a theorem of projective geometry.  The projective formulation is
Theorem \ref{projective lift theorem}.  There, orthogonality between the
edges of the reciprocal and facets of $M$ is replaced by duality
between lines and subspaces of codimension 1 in $\RP^{d}$.   Other
terminology is essentially as in the affine version.

\begin{theorem}
\label{projective lift theorem}
Let $\cM$ be a finite $(d-3)$-simple
homology $d$-manifold with $H_1(\cM,\Z_2)=0$. Let $M$ be a non-degenerate
generic piecewise-linear realization in $\RP^{d}$. Then $M$ has a sharp
piecewise-linear realization $L$ in $\RP^{d+1}$ that is projectively
equivalent to $M$. Such a realization is uniquely determined by choosing
$d+2$ free real parameters.
\end{theorem}

\begin{prop}In fact, as follows from recent work of Rybnikov (2003a,b), even if in
the above corollary $\cM$ is only $(d-3)$-simple, it is still the projection of a
convex polytope. This is a consequence of the fact that convexity at $(d-3)$-faces
(plus other van Heijenoort (1952) conditions ) is sufficient to establish global
 convexity.
\end{prop}
\section{Appendix}

\subsection{Counterexamples}

\begin{example} \label{X:infinite}{\rm
The graph is shown in Figure 1. We use the following notation for elements of
$Z_1(\G)$:
$$
\e_k = D_kL_kD_{k+1}U_{k}D_k, \quad \g_k = L_kL_{k+1}D_{k+1}L_k, \quad \h_k =
U_kD_{k+1}U_{k+1}U_k.
$$
Then $\{\e_k, \g_{k}, \h_{k} \}_{k \ge 1}$ is an integral circle basis of $\G$. As
usual, we denote an element of $Z_1(\G;\Z_2)$ corresponding to $\x \in Z_1(\G)$ by $[
\x ]$.  Let
$$
\mathbf{B}=\{\e_{k}+\g_{k}+\h_{k}\}_{k \ge 1} \cup \{\g_k+\e_{k+1} \}_{k \ge 1} \cup
\{\h_k+\e_{k+1}\}_{k \ge 1}
$$
and let $B=\{[\b] \mid \b \in \mathbf{B}\}$. $B$ is obviously a basis of
$Z_1(\G;\Z_2)$, because $[\e_k] =
[\e_{k}+\g_{k}+\h_{k}]+[\g_{k}+\e_{k+1}]+[\h_{k}+\e_{k+1}]$. Let $H$  denote the
subgroup of $Z_1(\G)$ (hence of $\Z E$)  generated by  $\mathbf{B}$, and let $\chi: \Z
E \to \Z E/H$ be the quotient map. Now, let us set the gain group to be $\Z E/H$ and
the gain map $g$ to be $\chi a$, where $a$ is the abelianization mapping from $F_E$  to
$\Z E$. Then the essential gain group $g(\pi(\G))=\chi(Z_1(\G))=Z_1(\G)/H$. By
construction, $B$ is a basis of the binary cycle space that consists of balanced
circles. If we show that $H \neq Z_1(\G)$,  then the constructed gain graph is not
balanced, because $Z_1(\G)/H=\langle \chi(\e_k), \chi(\g_{k}), \chi(\h_{k}) \mid k \ge
1 \rangle=g(\pi(\G))$.

\begin{figure}
\vbox to 3.3 truein{}
\begin{center}
\includegraphics{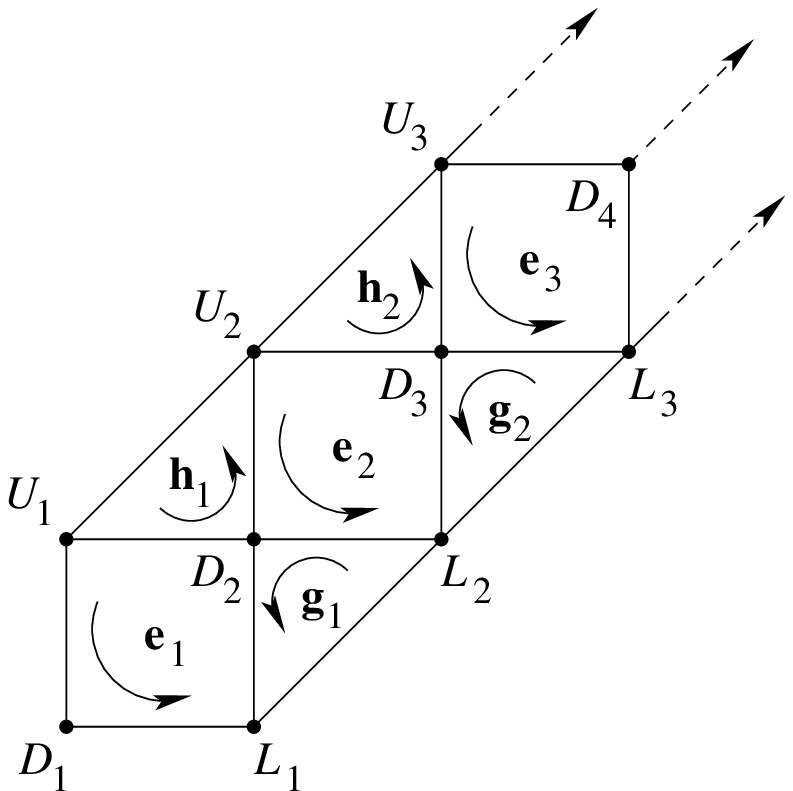}
\end{center}
\caption{Counterexample to the Circle Test with infinitely 2-divisible elements.}
\label{F:egh-graph:appendix}
\end{figure}

\begin{lemma}
In the above construction, $H \neq Z_1(\G)$.
\end{lemma}

\begin{proof}
The proof is based on the observation that there is no finite linear combination of
elements of $H$ that is equal to $\e_1$. Suppose there is such a combination. Then the
coefficient in front of $\e_{1}+\g_{1}+\h_{1}$ is 1. This, in turn, implies that the
coefficients in front of $\e_{2}+\h_{1}$ and $\e_{2}+\g_{1}$ are $-1$. Suppose we have
shown that the coefficient in front of $\e_{i}+\g_{i}+\h_{i}$ must be $(-2)^{i-1}$.
Then, obviously, the coefficients in front of $\e_{i+1}+\h_{i}$ and $\e_{i+1}+\g_{i}$
are $(-2)^{i-1}$. Therefore, the coefficient in front of $\e_{i+1}+\g_{i+2}+\h_{i+2}$
is $(-2)^{i}$. Thus, there is no finite linear combination of elements of $H$ that
gives $\e_1$.
\end{proof}

We want to prove that the essential gain group $g(\pi(\G))$ is 2-divisible. Denote by
$\Z [\frac{1}{2}]$ the ring of dyadic rationals, i.e., the ring of  rational numbers
whose denominators are powers of 2. Its additive group $\Z [\frac{1}{2}]^+$ is the
unique smallest group that contains an infinitely 2-divisible element of infinite
order.

\begin{lemma}
$Z_1(\G)/H$ is isomorphic to $\Z [\frac{1}{2}]^+$.
\end{lemma}

\begin{proof}
It is well known that $\Z [\frac{1}{2}]^{+}$ can be written as $\langle d_1,
d_2,\ldots\, |\, d_1=2d_2,\ d_2=2d_3,\ \ldots \rangle = A/D$, where $A$ is the free
abelian group with generators $d_1, d_2, \ldots$ and $D$ is its subgroup generated by
$d_i-2d_{i+1}$. Let $\phi$ be the map from $Z_1(\G)/H$ to $A/D$ defined by
$\phi(\e_i+H)=d_i+D$, $\phi(\g_i+H)=-d_{i+1}+D$, and $\phi(\h_i+H)=-d_{i+1}+D$ for $i
\ge 1$. It is easy to see that $\phi$ is an isomorphism.
\end{proof}
} \end{example}

Notice that a circle $\e_{i}+\g_{i}+\h_{i}$ has 5 edges. We do not know if one can
construct an example of an unbalanced infinite abelian gain graph with a torsion-free
gain group, where there is a basis of the binary cycle space that consists of balanced
circles with at most 4 edges. We suspect that such an example exists. However, we have
the following conjecture.

\begin{conjecture}
Let $G$ be a torsion-free abelian group, possibly 2-divisible, and suppose $Z_1(\G)$
has a basis that consists of balanced circles with at most 3  edges. Then $\GG$ is
balanced!
\end{conjecture}

There are also finite counterexamples, but the gain group necessarily has an element of
odd order.

\begin{figure}
\vbox to 5 truein{}
\begin{center}
\includegraphics{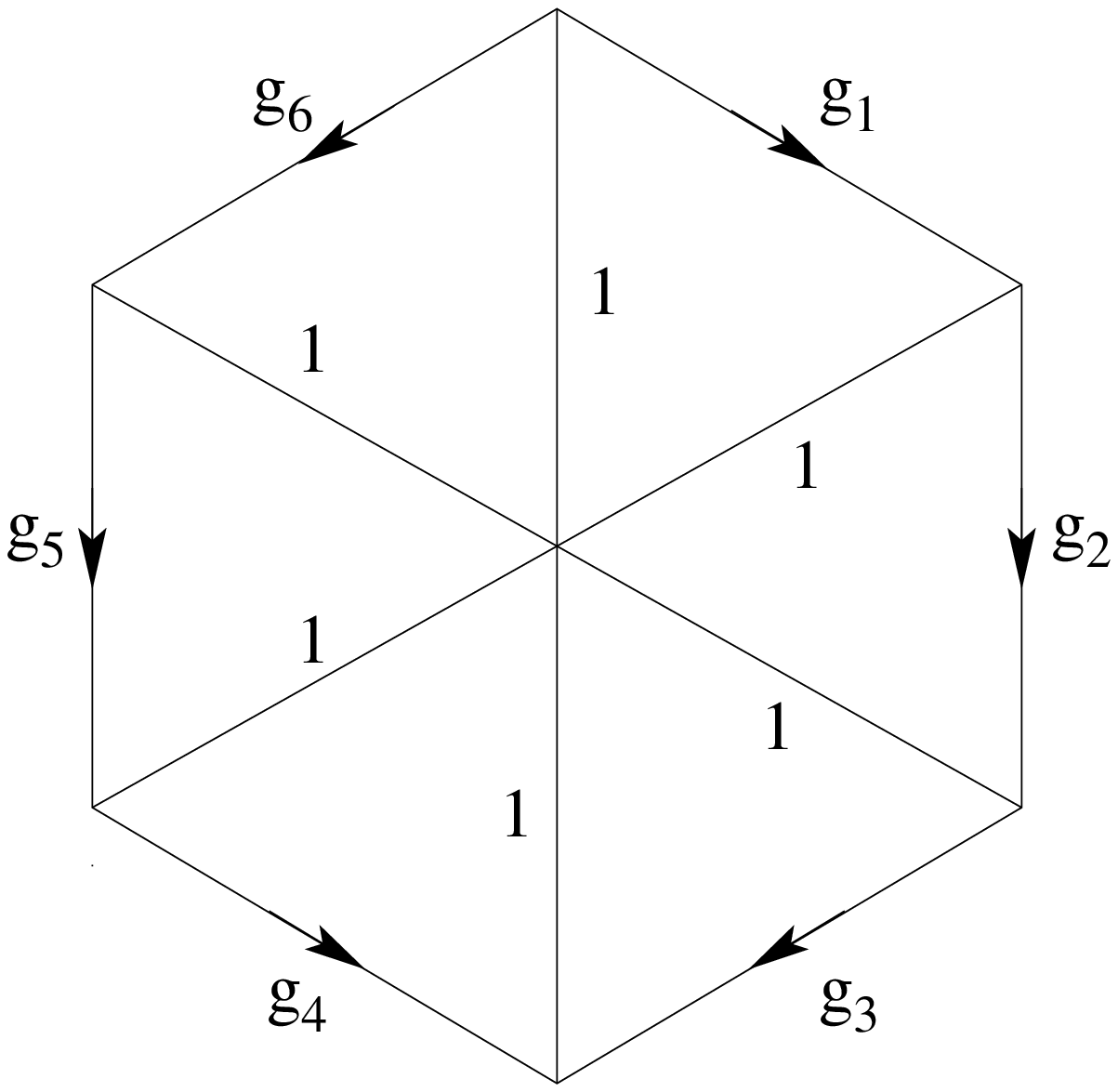}
\end{center}
\caption{The wheel $W_6$ with gains switched as in Example \ref{X:wheel}.  All products
$g_{i+1}g_{i+2}\cdots g_{i-1}=1$.}    \label{F:wheel}
\end{figure}

\begin{example} \label{X:wheel}
{\rm From Rybnikov and Zaslavsky (20xx), Theorem 6.16,
we take the example of the wheel $W_{2k}$ with gain group $\bbZ_{2k-1}$ (Figure
\ref{F:wheel}).  The binary cycle basis consists of the $2k$ Hamiltonian circles.
Switch so the spokes of the wheel have identity gain; then balance of the Hamiltonian
circles implies that any $2k-1$ consecutive edges along the outer circle have gain
product 1.  From this it follows that they all have the same gain $a$ and $a^{2k-1}=1$.
These gains with gain group $\bbZ_{2k-1}$ give an unbalanced gain graph that has a
binary cycle basis composed of balanced circles.
 }
\end{example}

Another example from Rybnikov and Zaslavsky (20xx), Theorem 6.16,
is $2C_{2k}$, an even circle with all edges doubled, with the same gain group.

\subsection{Joswig's Projectivities and Combinatorial Topology of Manifolds}
A special case of permutation gain graphs has been studied by Joswig (2002) and
Izmestiev \& Joswig (2003) in connection with coloring of tilings of manifolds and
problems in geometric topology of covering spaces.  Joswig considered gain graphs whose
underlying graphs appear as dual graphs of simplicial PL-manifolds. In this case, the
underlying graph $\G$ is the dual graph of some finite, pure, $d$-dimensional, strongly
connected simplicial  complex $\Delta$ on $n$ vertices, and the gain group is $S_n$,
the symmetric group on the vertices of $\Delta$. The gain map $g$ is determined by the
triangulation of $\Delta$ as follows: an edge $\sigma \tau$ of $\G$, where $\sigma$ and
$\tau$ are adjacent $d$-simplices of $\Delta$, is mapped to the transposition that
permutes two vertices of $\sigma \cup \tau$ that do not belong to their common facet.
Thus, $g(\pi(\G))$ is a subgroup of $S_n$, generated by a subset of transpositions on
$\{1,\ldots,n\}$ that correspond to adjacent $d$-simplices (such transpositions are
called ``perspectivities'' by Joswig). Let us fix a vertex $v$ of $\G$ (i.e. a
$d$-simplex $v=\{1,\ldots,d+1\}$ of $\Delta$) and consider the group of all closed
paths on $\G$ which  start and end at $v$; this group can be identified with $\pi(\G)$.
Then $g(\pi(\G))$ is the subgroup of all permutations that can be written as products
of transpositions along closed paths starting at $v$.  Such products of transpositions
are called ``projectivities'' by Joswig (2002) and their group is denoted by
$\Pi(\Delta,v)$ and referred to as the group of projectivities.  (Thus, Joswig's
$\Pi(\Delta,v)$ is the same as $g(\pi(\G))$ in our notation.) If $[v,v_1,\ldots,v_n,v]$
is a closed path in the dual graph starting from simplex $v$, then it is clear that
$g([v,v_1,\ldots,v_n,v_o])$, restricted to the vertices of $v$, is a permutation on
$\{1,\ldots,d+1\}$.   Therefore,  $g(\pi(\G))$ is a subgroup of the set-stabilizer
(also called isotropy subgroup) of ${1,\ldots,d+1}$. The main object of Joswig's (2002)
and Izmestiev \& Joswig's (2003) studies is the size and structure of this subgroup for
various combinatorial manifolds. This group was implicitly considered earlier, in the
special case where $\Delta$ is a 3-sphere, by Goodman and Onishi (1978).

\subsection{Remarks on Proof of Theorem 6.4}

Let us denote by $\textrm{C}^k (\G, \R^d)$ the $\R$-space of $k$-co-chains with values
in $\R^d$ of graph  $\G$. From the algebraic standpoint a reciprocal for $M$ is a
linear operator $r$ from the $\R$-space of 0-cochains on $\G(V,E)$ with values in
$\R^d$, with the property that $\delta r(C_i,C_j)$ is orthogonal to the oriented facet
$(C_i,C_j)$, where $\delta$ is the coboundary map from $ \textrm{C}^0(\G,\R^d)$ to
$\textrm{C}^1(\G,\R^d)$, defined by $\delta r(C_i,C_j) = r(C_j)-r(C_i)$ for each
oriented facet $(C_i,C_j)$ (i.e. an oriented edge of $\G$).

\end{document}